\documentclass[amstex,12pt,russian,amssymb]{article}
\usepackage{mathtext}
\usepackage[cp1251]{inputenc}
\usepackage[T2A]{fontenc}
\usepackage[russian]{babel}
\usepackage[dvips]{graphicx}
\usepackage{amsmath}
\usepackage{amssymb}
\usepackage{amsxtra}
\usepackage{latexsym}
\usepackage{ifthen}

\begin{document}

\newcounter{lemma}
\newcommand{\lemma}{\par \refstepcounter{lemma}%
{\bf Лемма \arabic{lemma}.}}

\newcounter{corollary}
\newcommand{\corollary}{\par \refstepcounter{corollary}%
{\bf Следствие \arabic{corollary}.}}

\newcounter{remark}
\newcommand{\remark}{\par \refstepcounter{remark}%
{\bf Замечание \arabic{remark}.}}

\newcounter{theorem}
\newcommand{\theorem}{\par \refstepcounter{theorem}%
{\bf Теорема \arabic{theorem}.}}

\newcounter{proposition}
\newcommand{\proposition}{\par \refstepcounter{proposition}%
{\bf Предложение \arabic{proposition}.}}

\newcounter{definition}
\newcommand{\definition}{\par \refstepcounter{definition}%
{\bf Определение \arabic{definition}.}}

\renewcommand{\refname}{\centerline{\bf Список литературы}}

\newcommand{\proof}{{\it Доказательство.\,\,}}

\noindent УДК 517.5

{\bf Д.С.~Доля} (Институт математики НАН Украины),

{\bf Е.А.~Севостьянов} (Житомирский государственный университет им.\
И.~Франко)

\medskip
{\bf Д.С.~Доля} (Інститут математики НАН України),

{\bf Є.О.~Севостьянов} (Житомирський державний університет ім.\
І.~Франко)

\medskip
{\bf D.S.~Dolya} (Institute of Mathematics of NAS of Ukraine),

{\bf E.A.~Sevost'yanov} (Zhitomir State University of I.~Franko)

\medskip
{\bf Об устранимых особенностях одного класса отображений,
удовлетворяющих модульным неравенствам}

{\bf Про усувні сингулярності одного класу відображень, що
задовольняють модульні нерівності}

{\bf On removable singularities of one class of mappings satisfying
moduli ine\-qua\-li\-ti\-es}

\medskip
В настоящей работе исследованы вопросы, связанные с локальным
поведением так называемых $Q$-отображений, включающих в себя классы
квазиконформных отображений и отображений с ограниченным искажением.
Показано, что, если отображение такого типа растёт в окрестности
изолированной точки границы не быстрее некоторой функции радиуса
шара, то эта точка является устранимой особой точкой отображения,
либо полюсом.

\medskip
В даній роботі досліджено питання, пов'язані з локальною поведінкою
так званих $Q$-відображень, що включають до себе класи
квазіконформних відображень та відображень з обмеженим спотворенням.
Доведено, що якщо відображення такого типу зростає в околі
ізольованої точки межі не швидше за деяку функцію радіусу кулі, то
ця точка є усувною сингулярністю відображення, або полюсом.

\medskip
A paper is devoted to study of local behavior of so-called
$Q$-mappings including qua\-si\-con\-for\-mal mappings and mappings
with bounded distortion. It is showed that, such mappings have
removable isolated singularities whenever the grow of the mappings
is note more than some function of a radius of a ball.

\newpage
{\bf 1. Введение.} Настоящая заметка посвящена поиску условий на
отображение $f:D\setminus\{b\}\rightarrow {\Bbb R}^n, $ заданное в
окрестности точки $b$ пространственной области $D\subset{\Bbb R}^n,$
$n\ge 2,$ при которых это отображение продолжается в точку $b$
непрерывным образом. Известно, что даже аналитическая функция
$\varphi:D\setminus\{b\}\rightarrow {\Bbb C},$ $D\subset {\Bbb C},$
заданная в области $D\setminus\{b\}$ комплексной плоскости ${\Bbb
C},$ вообще говоря, не продолжается в точку $b$ по непрерывности
($\varphi(z)=\exp\{1/z\},\quad b=0$). Однако, как было показано в
работе \cite{S$_5$}, такое продолжение имеет место для некоторого
класса отображений даже более общих, чем аналитические функции, в
том случае, если рост отображения в окрестности точки $b$ не слишком
велик. В случае отображений с ограниченным искажением см. также
работу \cite{Va} по этому поводу.

В настоящей работе класс отображений, для которых имеется подобное
непрерывное продолжение, будет ещё более расширен. Здесь мы покажем,
что для утверждения такого рода достаточно обобщённого неравенства
типа Полецкого, определяющего открытые дискретные кольцевые
$Q$-отображения (см. \cite{MRSY$_1$}--\cite{MRSY}, \cite{BGMV},
\cite{GRSY}). При этом, на порядок роста самих отображений, к
сожалению, необходимо потребовать несколько более жёсткие условия
аналогичного же характера (здесь имеются в виду
степенно-логарифмические ограничения, где <<порядок степени>> может
быть заранее неизвестен). Всё сказанное, что следует отдельно
отметить, не является содержательным в случае гомеоморфизмов (см.,
напр., \cite{IR}). В этом случае, никаких ограничений на рост
отображений не требуется, кроме, разумеется, аналитических условий
на функцию $Q,$ отвечающую за меру искажения модуля семейств кривых
при отображении.

\medskip
Теперь более подробно о содержательной части работы и её основных
результатах. Всюду далее $D$ -- область в ${\Bbb R}^n,$ $n\ge 2,$
$m$ -- мера Лебега ${\Bbb R}^n,$ ${\rm dist\,}(A,B)$ --
ев\-кли\-дово расстояние между множествами $A, B\subset {\Bbb R}^n,$
${\rm dist\,}(A,B)=\inf\limits_{x\in A, y\in B} |x-y|,$ $(x,y)$
обозначает (стандартное) скалярное произведение векторов $x,y\in
{\Bbb R}^n,$ ${\rm diam\,}A$ -- ев\-кли\-дов диаметр множества
$A\subset {\Bbb R}^n,$ $B(x_0, r)=\left\{x\in{\Bbb R}^n: |x-x_0|<
r\right\},$ ${\Bbb B}^n := B(0, 1),$ $S(x_0,r) = \{ x\,\in\,{\Bbb
R}^n : |x-x_0|=r\},$ ${\Bbb S}^{n-1}:=S(0, 1),$ $\omega_{n-1}$
означает площадь сферы ${\Bbb S}^{n-1}$ в ${\Bbb R}^n,$ $\Omega_n$
-- объём единичного шара ${\Bbb B}^{n}$ в ${\Bbb R}^n,$ запись
$f:D\rightarrow {\Bbb R}^n$ предполагает, что отображение $f,$
заданное в области $D,$ непрерывно. Отображение $f:D\rightarrow
{\Bbb R}^n$ называется {\it дискретным}, если прообраз
$f^{-1}\left(y\right)$ каждой точки $y\,\in\,{\Bbb R}^n$ состоит из
изолированных точек, и {\it открытым}, если образ любого открытого
множества $U\subset D$ является открытым множеством в ${\Bbb R}^n.$

\medskip {\it Кривой} $\gamma$ мы называем непрерывное
отображение отрезка $[a,b]$ (открытого интервала $(a,b),$ либо
полуоткрытого интервала вида  $[a,b)$ или $(a,b]$) в ${\Bbb R}^n,$
$\gamma:[a,b]\rightarrow {\Bbb R}^n.$ Под семейством кривых $\Gamma$
подразумевается некоторый фиксированный набор кривых $\gamma,$ а
$f(\Gamma)=\left\{f\circ\gamma|\gamma\in\Gamma\right\}.$ Следующие
определения могут быть найдены, напр., \cite[разд.~1--6, гл.~
I]{Va$_1$}. Борелева функция $\rho:{\Bbb R}^n\,\rightarrow
[0,\infty]$ называется {\it допустимой} для семейства $\Gamma$
кривых $\gamma$ в ${\Bbb R}^n,$ если криволинейный интеграл первого
рода от функции $\rho$ по каждой (локально спрямляемой) кривой
$\gamma\in \Gamma$ удовлетворяет условию $\int\limits_{\gamma}\rho
(x)|dx|\ge 1.$ В этом случае мы пишем: $\rho \in {\rm adm}\Gamma.$
{\it Модулем} семейства кривых $\Gamma $ называется величина
$M(\Gamma)=\inf_{\rho \in \,{\rm adm}\,\Gamma} \int\limits_D \rho ^n
(x)\ \ dm(x).$ Свойства модуля в некоторой мере аналогичны свойствам
меры Лебега $m$ в ${\Bbb R}^n.$ Именно, модуль пустого семейства
кривых равен нулю, $M(\varnothing)=0,$ модуль обладает свойством
монотонности относительно семейств кривых $\Gamma_1$ и $\Gamma_2:$
$\Gamma_1\subset\Gamma_2\Rightarrow M(\Gamma_1)\le M(\Gamma_2),$ а
также свойством полуаддитивности:
$M\left(\bigcup\limits_{i=1}^{\infty}\Gamma_i\right)\le
\sum\limits_{i=1}^{\infty}M(\Gamma_i),$ см.
\cite[теорема~6.2]{Va$_1$}. Говорят, что некоторое свойство
выполнено для {\it почти всех (п.в.) кривых} области $D$, если оно
имеет место для всех кривых, лежащих в $D$, кроме (возможно)
некоторого их подсемейства, модуль которого равен нулю. Далее символ
$\Gamma(E,F,D)$ означает семейство всех кривых
$\gamma:[a,b]\rightarrow{\Bbb R}^n,$ которые соединяют $E$ и $F$ в
$D,$ т.е. $\gamma(a)\in E,$ $\gamma(b)\in F$ и $\gamma(t)\in D$ при
$t\in (a, b).$

\medskip
Пусть $r_0\,=\,{\rm dist}\, (x_0\,,\partial D),$
$Q:D\rightarrow\,[0\,,\infty]$ -- измеримая по Лебегу функция,
$S_{\,i}\,=\,S(x_0,r_i),$ $i=1,2.$ Отображение
$f:D\setminus\{x_0\}\rightarrow \overline{{\Bbb R}^n},$
$\overline{{\Bbb R}^n}:={\Bbb R}^n\cup\{\infty\},$ называется {\it
кольцевым $Q$-отоб\-ра\-же\-нием в точке $x_0\,\in\,D,$} если
соотношение
$$M\left(f\left(\Gamma\left(S_1,\,S_2,\,A\right)\right)\right)\ \le
\int\limits_{A} Q(x)\cdot \eta^n(|x-x_0|)\ dm(x)$$
выполнено для любого кольца $A=A(r_1,r_2, x_0)=\{x\in {\Bbb R}^n:
r_1<|x-x_0|<r_2\},$  $0<r_1<r_2< r_0,$ и для каждой измеримой
функции $\eta : (r_1,r_2)\rightarrow [0,\infty ]\,$ такой, что
\begin{equation}\label{eq*3}
\int\limits_{r_1}^{r_2}\eta(r)\ dr\ \ge\ 1\,.
\end{equation}

\medskip
Всюду далее $q_{x_0}(r)$ означает среднее интегральное значение
$Q(x)$ над сферой $S(x_0, r),$
$$
q_{x_0}(r):=\frac{1}{\omega_{n-1}r^{n-1}}\int\limits_{S(x_0, r)}
Q(x)\,dS\,,
$$
где $dS$ -- элемент площади поверхности $S.$

\medskip Напомним, что
изолированная точка $x_0$ границы $\partial D$ области $D$
называется {\it устранимой} для отображения $f,$ если существует
конечный предел $\lim\limits_{x\rightarrow x_0}\,f(x).$ Если
$f(x)\rightarrow \infty$ при $x\rightarrow x_0,$ точку $x_0$ будем
называть {\it полюсом.} Изолированная точка $x_0$ границы $\partial
D$ называется {\it существенно особой точкой} отображения
$f:D\rightarrow {\Bbb R}^n,$ если при $x\rightarrow x_0$ нет ни
конечного, ни бесконечного предела. Основной результат настоящей
статьи заключает в себе следующее

\medskip
\begin{theorem}\label{th1A} {\sl Пусть $Q:D\rightarrow\,(0\,,\infty]$ -- измеримая по Лебегу функция,
$b\in D$ и $f:D\setminus\{b\}\rightarrow {\Bbb R}^n$ --
открытое дискретное кольцевое $Q$-отображение в точке $b.$
Зафиксируем $\varepsilon_0>0,$ такое что $\varepsilon_0<{\rm
dist}\,(b,
\partial D).$ Предположим, что
$\int\limits_{0}^{\varepsilon_0}\frac{dt}{tq_{b}^{1/(n-1)}(t)}=\infty,$
$q_{b}(t)<\infty$ при почти всех $t\in (0, \varepsilon_0).$ Тогда,
если при $0<p<\left(\frac{1}{4}\right)^{1/(n-1)}$ и некоторой
постоянной $C>0$
\begin{equation}\label{eq2}
|f(x)|\le C \cdot
\exp\left\{p\int\limits_{|x-b|}^{\varepsilon_0}\frac{dt}{tq_{b}^{1/(n-1)}(t)}\right\}\qquad\forall\,\,
x\in B(b, \varepsilon_0)\setminus \{b\}\,,
\end{equation}
то отображение $f$ имеет предел (конечный либо бесконечный) в точке
$b.$ Если же вместо условия (\ref{eq2}) имеет место более сильное
предположение
$$\lim\limits_{x\rightarrow b} |f(x)|\cdot
\exp\left\{-p\int\limits_{|x-b|}^{\varepsilon_0}\frac{dt}{tq_{b}^{1/(n-1)}(t)}\right\}
=0\,,\quad 0<p< \left(\frac{1}{4}\right)^{1/(n-1)}\,,$$
то точка $x=b$ является для отображения $f$ устранимой особой
точкой, исключая полюс.}
\end{theorem}

\medskip
Отметим, что для отображений с ограниченным искажением (когда
функция $Q(x)$ ограничена) теорема \ref{th1A} была получена в работе
\cite[теорема~4.2]{Va}.

\medskip
{\bf 2. Основная лемма.}  Компактное множество $G\subset
\overline{{\Bbb R}^n}$ условимся называть {\it множеством нулевой
ёмкости,} пишем ${\rm cap\,}G =0,$ если существует компакт $T\subset
\overline{{\Bbb R}^n}\setminus G,$ $\overline{\overline{{\Bbb
R}^n}\setminus T}\ne \overline{{\Bbb R}^n},$ такой что $M(\Gamma(T,
G, {\Bbb R}^n))=0,$ см., напр., \cite[разд.~2, гл.~III и
предложение~10.2, гл.~II]{Ri}. Будем говорить, что произвольное
множество $G$ имеет ёмкость нуль, если произвольное его компактное
подмножество $G_0$ имеет нулевую ёмкость. Множества ёмкости нуль,
как известно, всюду разрывны (любая компонента их связности
вырождается в точку), т.е., условие ${\rm cap\,}G =0$ влечёт, что
${\rm Int\,}G=\varnothing,$ см., напр., \cite[следствие~2.5, гл.~
III]{Ri}. Открытое множество $U\subset D,$ $\overline{U}\subset D,$
называется {\it нормальной окрестностью} точки $x\in D$ при
отображении $f:D\rightarrow {\Bbb R}^n,$ если $U\cap
f^{\,-1}\left(f(x)\right)=\left\{x\right\}$ и $\partial
f(U)=f(\partial U),$ см., напр., \cite[разд.~4, гл.~I]{Ri}.

\medskip
\begin{proposition}\label{pr2}{\sl\,
Пусть $f:D\rightarrow {\Bbb R}^n$ открытое дискретное отображение.
Тогда для каждого $x\in D$ существует $s_x,$ такое, что при всех
$s\in (0, s_x)$ компонента связности множества $f^{-1}\left(B(f(x),
s)\right),$ содержащая точку $x,$ и обозначаемая символом
$U(x,f,s),$ является нормальной окрестностью точки $x$ при
отображении $f,$ при этом $f(U(x,f,s))=B(f(x), s)$ и ${\rm
diam\,}U(x,f,s)\rightarrow 0$ при $s\rightarrow 0,$ см., напр.,
лемму 4.9 гл. I в \cite{Ri}.}
\end{proposition}

\medskip
Для отображения $f:D\rightarrow{\Bbb R}^n,$ множества $E\subset D$ и
$y\in{\Bbb R}^n,$  определим {\it функцию кратности $N(y,f,E)$} как
число прообразов точки $y$ во множестве $E,$ т.е.
$N(y,f,E)={\rm card}\,\left\{x\in E: f(x)=y\right\}.$ Для
доказательства основных результатов работы нам необходимо
воспользоваться следующим утверждением, см. \cite[лемма~5.1]{SS}.

\medskip
\begin{proposition}\label{pr4}{\sl Пусть
$Q:{\Bbb B}^n\rightarrow [0, \infty]$ -- измеримая по Лебегу
функция, $f:{\Bbb B}^n\setminus \left\{0\right\} \rightarrow
\overline{{\Bbb R}^n},$ $n \ge 2\,,$ -- открытое дискретное
кольцевое $Q$-отоб\-ра\-же\-ние. Пусть, кроме того,
$${\rm cap}\,\left(\overline{{\Bbb R}^n}\setminus\,f\left({\Bbb
B}^n\setminus\left\{0\right\}\right)\right)>0\,.$$ Предположим, что
существует $\varepsilon_0\in(0,1)$ такое, что при
$\varepsilon\rightarrow 0$
$$
\int\limits_{\varepsilon<|x|<\varepsilon_0}Q(x)\cdot\psi^n(|x|) \
dm(x)\,=\,o\left(I^n(\varepsilon, \varepsilon_0)\right)\,,
$$
где $\psi(t)$ -- измеримая по Лебегу функция, такая что $\psi(t)>0$
п.в., и
$$
I(\varepsilon,
\varepsilon_0)=\int\limits_{\varepsilon}^{\varepsilon_0}\psi(t)dt <
\infty
$$
для всех $\varepsilon\in(0, \varepsilon_0).$ Тогда $f$ имеет
непрерывное продолжение $f:{\Bbb B}^n\rightarrow\overline{{\Bbb
R}^n}$ в ${\Bbb B}^n.$ }
\end{proposition}

\medskip
Важнейшую роль при доказательстве основных результатов работы играет
следующая лемма.

\medskip
\begin{lemma}\label{lem1}{\sl Пусть $Q:D\rightarrow\,(0\,,\infty]$ -- измеримая по Лебегу функция,
$b\in D$ и $f:D\setminus\{b\}\rightarrow {\Bbb R}^n$ -- открытое
дискретное кольцевое $Q$-отображение в точке $b.$ Зафиксируем
$\varepsilon_0>0,$ такое что $\varepsilon_0<{\rm dist}\,(b,
\partial D).$ Предположим, что для некоторых постоянной $A>0,$
измеримой по Лебегу функции $\psi:(0, \varepsilon_0)\rightarrow (0,
\infty)$ и $\varepsilon\rightarrow 0$
\begin{equation} \label{eq4}
\int\limits_{\varepsilon<|x-b|<\varepsilon_0}Q(x)\cdot\psi^n(|x-b|)
\ dm(x)\le \frac{A\cdot I^n(\varepsilon,
\varepsilon_0)}{\left(\log\varphi(\varepsilon)\right)^{n-1}}\,,
\end{equation} где
\begin{equation} \label{eq11}
I(\varepsilon, \varepsilon_0):
=\int\limits_{\varepsilon}^{\varepsilon_0}\psi(t)dt < \infty \qquad
\forall\quad\varepsilon \in(0, \varepsilon_0)\,.\end{equation}
Зафиксируем строго убывающую функцию $\varphi:(0, \infty)\rightarrow
(0, \infty),$ для которой $\varphi(t)\rightarrow \infty$ при
$t\rightarrow 0.$
Тогда, если при некотором $0<p<
\left(\frac{\omega_{n-1}}{4A}\right)^{1/(n-1)}$ и некоторой
постоянной $C>0$ выполняется условие
\begin{equation}\label{eq3}
|f(x)|\le C \cdot \varphi^p (|x-b|)\quad \forall\,x\in
B(0,\varepsilon_0)\setminus\{0\}\,,
\end{equation}
то отображение $f$ имеет предел (конечный либо бесконечный) в точке
$b.$ }
\end{lemma}

\begin{proof} Предположим противное, а именно, что точка $b$
является существенно особой точкой отображения $f.$ Не ограничивая
общности рассуждений, можно считать, что $b=0$ и $C=1.$ Поскольку
сфера $S(0, \varepsilon_0)$ является компактным множеством в
$D\setminus\{0\},$ найдётся $R>0,$ такое что
\begin{equation}\label{eq10}
f\left(S(0, \varepsilon_0)\right)\subset B(0, R)\,.
\end{equation}
Поскольку $b=0$ является существенно особой точкой отображения $f,$
в виду условия (\ref{eq4}) и предложения \ref{pr4} отображение $f$ в
$B(0, \varepsilon_0)\setminus \{0\}$ принимает все значения в ${\Bbb
R}^n,$ за исключением, может быть, некоторого множества ёмкости
нуль, т.е., $N\left(y, f, B(0, \varepsilon_0)\setminus
\{0\}\right)=\infty$ при всех $y\in {\Bbb R}^n\setminus E,$ где
${\rm cap\,}E=0.$ Так как $E$ имеет ёмкость нуль, множество ${\Bbb
R}^n\setminus E$ не может быть ограниченным. В таком случае,
найдётся $y_0\in {\Bbb R}^n\setminus \left(E\cup B(0, R)\right)$ а,
значит, найдётся точка $x_1\in f^{-1}(y_0)\cap (B(0,
\varepsilon_0)\setminus\{0\}).$ По предложению \ref{pr2}, при
некотором фиксированном $r>0,$ такая точка имеет нормальную
окрестность $U_1:=U(x_1, f , r),$ такую, что $f(U_1)=B(y_0, r).$

Полагаем $d:={\rm dist\,}\left(0, \overline{U_1}\right).$ Пусть
$a\in (0, d)$ и $V:=B(0, \varepsilon_0)\setminus\overline{B(0, a)}.$
В силу неравенства (\ref{eq3}), строгого убывания функции $\varphi,$
а также предположения о том, что $C=1,$ имеем
\begin{equation}\label{eq6}
f(V)\subset B\left(0, \varphi^p(a)\right)\,.
\end{equation}
Поскольку $z_0:=y_0+re\in \overline{B(y_0,
r)}=f\left(\overline{U(x_1, f, r)}\right),$ будем иметь: $z_0\in
f(V).$ Следовательно, найдётся последовательность точка
$\widetilde{x_1}\in \overline{U_1},$ такая что
$f(\widetilde{x_j})=z_0.$ Обозначим через $H$ полусферу
$H=\left\{e\in {\Bbb S}^{n-1}: (e, y_0)>0\right\},$ через
$\Gamma^{\,\prime}$ -- семейство всех кривых $\beta:\left[r,
\varphi^p(a)\right)\rightarrow {\Bbb R}^n$ вида $\beta(t)=y_0+te,$
$e\in H,$ а через $\Gamma$ -- семейство максимальных под\-ня\-тий
$\Gamma^{\,\prime}$ при отображении $f$ относительно области $V$ с
началом в точке $\widetilde{x_1},$ которые существуют в силу
\cite[теорема~3.2, гл.~II]{Ri}.

\medskip
При любом фиксированном $e\in H,$ покажем, что для каждой кривой
$\beta=y_0+te$ и каждого максимального её поднятия $\alpha(t):[r,
c)\rightarrow V$ с началом в точке $\widetilde{x_1},$ $\alpha\in
\Gamma,$  существует последовательность $r_k\in [r, c),$ такая что
$r_k\rightarrow c-0$ при $k\rightarrow \infty$ и ${\rm
dist\,}(\alpha(r_k),
\partial V)\rightarrow 0$ при $k\rightarrow \infty.$
Предположим противное, тогда найдётся $e_0\in H,$ такое, что кривая
$\alpha(t),$ $t\in [r, c),$ являющаяся максимальным поднятием кривой
$\beta=y_0+te_0,$ лежит внутри $V$ вместе со своим замыканием. Пусть
$C(c,\,\alpha(t))$ обозначает предельное множество кривой $\alpha$
при $t\rightarrow c-0,$ тогда для каждого $x\in C(c,\,\alpha(t))$
найдётся последовательность $t_k\rightarrow \infty,$ такая, что
$x=\lim\limits_{k\rightarrow \infty}\alpha(t_k).$ По непрерывности
$f,$ поскольку, по предположению, $C(c,\,\alpha(t))\subset V,$ будем
иметь $f(x)=f(\lim\limits_{k\rightarrow \infty}\alpha(t_k))=
\lim\limits_{k\rightarrow \infty}\beta(t_k)= \beta(c),$ откуда
следует, что отображение $f$ постоянно на множестве
$C(c,\,\alpha(t)).$ Так как по условию $f$ дискретно, а множество
$C(c,\,\alpha(t)),$ очевидно, является связным, будем иметь
$C(c,\,\alpha(t))=p_1\in V.$

Полагаем $b_0:=\varphi^p(a).$ Пусть $c\ne b_0.$ Тогда можно
построить новое максимальное поднятие $\alpha^{\,\prime}$ кривой
$\beta$ с началом в точке $p_1.$ Объединяя поднятия $\alpha$ и
$\alpha^{\,\prime},$ получаем ещё одно поднятие
$\alpha^{\,\prime\prime}$ кривой $\beta$ с началом в точке
$\widetilde{x_{j_0}},$ что противоречит свойству максимальности
исходного поднятия $\alpha.$ Значит, $c=b_0.$

В таком случае, $C(b_0, \alpha(t))$ континуум внутри $V,$ при этом,
$C(b_0, \alpha(t))=p_1^{\,\prime}\in V$ и, значит, $\alpha$
продолжается до замкнутой кривой, определённой на отрезке $\left[r,
\varphi^p(a)\right].$ Обозначим эту кривую снова через $\alpha$
(обозначения не меняем). Тогда при всех $t\in \left[r,
\varphi^p(a)\right]$ имеем $\beta(t)=f(\alpha(t))\subset f(V),$ в
частности, полагая $t:=\varphi^p(a),$ рассмотрим элемент $z_1,$
определяемый по правилу $z_1:=y_0+\varphi^p(a)e_0.$  Ввиду включения
(\ref{eq6}), имеем
\begin{equation}\label{eq7}
z_1=y_0+\varphi^p(a)e_0\in f(V)\subset B\left(0,
\varphi^p(a)\right)\,.
\end{equation}
Однако, поскольку $e_0\in H,$
$$|z_1|=\left|y_0+\varphi^p(a)e_0\right|=\sqrt{|y_0|^2 +
2\left(y_0, \varphi^p(a)e_0\right)+\varphi^{2p}(a)}\ge$$
\begin{equation}\label{eq8}
\ge \sqrt{|y_0|^2 + \varphi^{2p}(a)}\ge \varphi^p(a)\,.
\end{equation}
Однако, соотношение (\ref{eq8}) противоречит (\ref{eq7}), что, в
свою очередь, опровергает предположение о включении замыкания кривой
$\alpha(t)$ во множество $V.$

Следовательно, ${\rm dist\,}(\alpha(r_k),
\partial V)\rightarrow 0$ при $k\rightarrow c-0$ и некоторой
последовательности $r_k\in [r, c),$ такой что $r_k\rightarrow c-0$ и
$k\rightarrow \infty,$ что и требовалось установить.

Заметим, что ситуация, когда ${\rm dist\,}(\alpha(r_k), S(0,
\varepsilon_0))\rightarrow 0$ при $k\rightarrow c-0,$ исключена.
Действительно, пусть эта ситуация имеет место. Тогда найдутся
$p_2\in S(0, \varepsilon_0)$ и подпоследовательность номеров $k_l,$
$l\in {\Bbb N},$ такие, что $\alpha(r_{k_l})\rightarrow p_2$ при
$l\rightarrow \infty.$ Отсюда, по непрерывности $f,$ получаем, что
$\beta(r_{k_l})\rightarrow f(p_2)$ при $l\rightarrow \infty,$ что
невозможно ввиду соотношения (\ref{eq10}), поскольку, при каждом
фиксированном $e\in H$ и $t\in \left[r, \varphi^p(a)\right),$ имеем
$|\beta(t)|=|y_0+te|=\sqrt{|y_0|^2 + 2t(y_0, e)+t^2}\ge |y_0|> R$ по
выбору $y_0.$

Из сказанного выше следует, что найдётся последовательность $r_k\in
[r, c),$ такая что $r_k\rightarrow c-0$ при $k\rightarrow \infty,$ и
$\alpha(r_k)\rightarrow p_3\in S(0, a).$ Кроме того, каждая такая
кривая $\alpha\in \Gamma$ пересекает сферу $S(0, d),$ поскольку,
согласно построению, $\alpha$ имеет начало вне шара $B(0, d).$ В
силу сказанного, $\Gamma>\Gamma(S(0, a+\varepsilon), S(0, d),
A(a+\varepsilon, d, 0))$ при сколь угодно малых $\varepsilon>0.$
Заметим, что $\Gamma^{\,\prime}>f(\Gamma),$ поэтому по определению
кольцевого $Q$-отображения в нуле
$$M(\Gamma^{\,\prime})\le M(f(\Gamma))\le M(\Gamma(S(0, a+\varepsilon), S(0, d),
A(a+\varepsilon, d, 0)))\le$$
\begin{equation}\label{eq9} \le \int\limits_{A(a+\varepsilon, d, 0)} Q(x)\cdot
\eta^n (x) dm(x)
\end{equation}
для каждой измеримой по Лебегу функции $\eta,$ удовлетворяющей
соотношению (\ref{eq*3}) при $r_1:=a+\varepsilon$ и $r_2:=d_0.$

\medskip
Рассмотрим функцию
$$\eta_{a+\varepsilon}(t)= \left\{
\begin{array}{rr}
\psi(t)/I(a+\varepsilon, d), &   t\in (a+\varepsilon, d),\\
0,  &  t\in {\Bbb R} \setminus (a+\varepsilon, d)
\end{array}
\right.\,, $$
где величина $I(a+\varepsilon, d)$ определена также, как в
(\ref{eq11}), а $\psi$ -- функция из условия леммы. Заметим, что
функция $\eta_{a+\varepsilon}(t)$ является измеримой по Лебегу,
кроме того,
$$\int\limits_{a+\varepsilon}^d \eta_{a+\varepsilon}(t) dt=1\,.$$ В таком случае, из
соотношения (\ref{eq9}) получаем:
$$M(\Gamma^{\,\prime} )\quad\le\quad \frac{1}{I^n(a+\varepsilon,
d)}\quad\int\limits_{a+\varepsilon<|x|<d} Q(x)\cdot \psi^n
(|x|)dm(x)\le$$
$$\le \frac{I^n(a+\varepsilon, \varepsilon_0)}{I^n(a+\varepsilon,
d)\cdot I^n(a,
\varepsilon_0)}\quad\int\limits_{a+\varepsilon<|x|<\varepsilon_0}
Q(x)\cdot \psi^n (|x|)dm(x)=$$$$=\left(1+\frac{I(d,
\varepsilon_0)}{I(a+\varepsilon, d)}\right)^n\frac{1}{\cdot
I^n(a+\varepsilon,
\varepsilon_0)}\quad\int\limits_{a+\varepsilon<|x|<\varepsilon_0}
Q(x)\cdot \psi^n (|x|)dm(x)\le$$
\begin{equation}\label{eq21}\le \frac{2}{I^n(a+\varepsilon,
\varepsilon_0)}\quad\int\limits_{a+\varepsilon<|x|<\varepsilon_0}
Q(x)\cdot \psi^n (|x|)dm(x)\end{equation} при всех $a+\varepsilon\in
(0, d_1)$ и некотором $d_1,$ $d_1\le d,$ поскольку, в силу
соотношения (\ref{eq4}), $I^n(a+\varepsilon,d)\rightarrow\infty$ при
$a+\varepsilon\rightarrow \infty.$ Снова, из (\ref{eq4}) и
(\ref{eq21}) получаем, что при $a+\varepsilon\in (0, d_1)$
\begin{equation}\label{eq12}
M(\Gamma^{\,\prime} )\le
\frac{2A}{\left(\log\varphi(a+\varepsilon)\right)^{n-1}}\,.
\end{equation}
С другой стороны, в силу \cite[разд.~7.7]{Va$_1$},
\begin{equation}\label{eq13}
M(\Gamma^{\prime})=\frac{1}{2}\frac{\omega_{n-1}}
{\left(\log\frac{\varphi^p(a+\varepsilon)}{r}\right)^{n-1}}\,.
\end{equation}
Тогда из неравенств (\ref{eq12}) и (\ref{eq13}) получаем:
$$\frac{1}{2}\frac{\omega_{n-1}}
{\left(\log\frac{\varphi^p(a+\varepsilon)}{r}\right)^{n-1}}\le
\frac{2A}{\left(\log\varphi(a+\varepsilon)\right)^{n-1}},$$
откуда
$$\left(\log\left(\frac{\varphi^p(a+\varepsilon)}{r}\right)^{
\left(\frac{2}{\omega_{n-1}}\right)^{\frac{1}{n-1}}}\right)^{n-1}\ge
\left(\log\left(\varphi(a+\varepsilon)\right)^{\left(\frac{1}{2A}\right)^
{\frac{1}{n-1}}}\right)^{n-1}\,,
$$
$$\log\left(\frac{\varphi^p(a+\varepsilon)}{r}\right)^{
\left(\frac{2}{\omega_{n-1}}\right)^{\frac{1}{n-1}}}\ge
\log\left(\varphi(a+\varepsilon)\right)^{\left(\frac{1}{2A}\right)^
{\frac{1}{n-1}}}\,,
$$
$$\frac{1}{r^{{
\left(\frac{2}{\omega_{n-1}}\right)^{\frac{1}{n-1}}}}}\ge
\left(\varphi(a+\varepsilon)\right)^{{\left(\frac{1}{2A}\right)^{\frac{1}{n-1}}}-p{
\left(\frac{2}{\omega_{n-1}}\right)^{\frac{1}{n-1}}}}\,.$$
Поскольку по выбору $0<p<
\left(\frac{\omega_{n-1}}{4A}\right)^{1/(n-1)},$ в правой части
последнего соотношения величина $\varphi(a+\varepsilon)$ берётся в
некоторой положительной степени. Переходя здесь к пределу при
$a+\varepsilon\rightarrow 0$ и учитывая, что по условию леммы
$\varphi(a)\rightarrow\infty$ при $a\rightarrow 0,$ получаем, что
$$\frac{1}{r^{{
\left(\frac{2}{\omega_{n-1}}\right)^{\frac{1}{n-1}}}}}\ge
\infty\,,$$ что невозможно.
Полученное противоречие означает, что точка $b=0$ не может быть
существенно особой для отображения $f.$ $\Box$
\end{proof}

\medskip
Отдельный случай представляет собой ситуация, когда $I(\varepsilon,
\varepsilon_0)\le M\cdot \log\varphi(\varepsilon)$ при некоторой
постоянной $M>0$ и $\varepsilon\rightarrow 0.$ Следующее утверждение
может быть получено из \cite[лемма~5]{S$_3$}.

\medskip
\begin{proposition}\label{pr3}
{\sl  Предположим, что $b\in D,$ $f:D\rightarrow B(0, R)$ --
открытое дискретное кольцевое $Q$-отображение в точке $b,$ при этом,
существуют измеримая по Лебегу функция $Q:D\rightarrow [0, \infty],$
числа $\varepsilon_0>0,$ $\varepsilon_0<{\rm dist\,}\left(b,
\partial D\right),$ и $A>0,$ такие, что
при $\varepsilon\rightarrow 0$ имеют место соотношения
(\ref{eq4})--(\ref{eq11}). Пусть, кроме того, существует постоянная
$M>0$ и $\varepsilon_1>0,$ $\varepsilon_1\in (0, \varepsilon_0),$
такие что при всех $\varepsilon\in (0, \varepsilon_1)$ выполнено
условие
\begin{equation} \label{eq4A}
I(\varepsilon, \varepsilon_0)\le M\cdot\log\varphi(\varepsilon)\,,
\end{equation}
где $I(\varepsilon, \varepsilon_0)$ определяется соотношением
(\ref{eq11}), а $\varphi:(0, \infty)\rightarrow [0, \infty)$ --
некоторая функция. Тогда при всех $x\in B(b, \varepsilon_1)$ имеет
место оценка
\begin{equation}\label{eq15}
|f(x)-f(b)|\le
\frac{\alpha_n(1+R^2)}{\delta}\exp\{-\widetilde{\beta_n}
I\left(|x-b|, \varepsilon_0\right)\}\,,
\end{equation}
где постоянные $\alpha_n$ и
$\widetilde{\beta_n}=\left(\frac{\omega_{n-1}}{AM^{n-1}}\right)^{1/(n-1)}$
зависят только от $n,$ а $\delta$ -- от $R.$ }
\end{proposition}

\medskip
\begin{corollary}\label{cor1}{ Предположим, что в условиях леммы \ref{lem1},
помимо соотношений (\ref{eq4}) и (\ref{eq11}) имеет место условие
(\ref{eq4A}), а вместо условия (\ref{eq3}) имеет место более сильное
предположение:
\begin{equation}\label{eq16}
\lim\limits_{x\rightarrow b}|f(x)|\cdot \exp\{-\beta_n I\left(|x-b|,
\varepsilon_0\right)\}=0\,,
\end{equation}
где $\beta_n=\left(\frac{\omega_{n-1}}{qAM^{n-1}}\right)^{1/(n-1)},$
$q>4$ и $M\ge 1.$ Тогда точка $x=b$ является устранимой для
отображения $f.$ }
\end{corollary}

\begin{proof}
Не ограничивая общности рассуждений, можно считать, что $b=0.$
Заметим, что
$\beta_n=\left(\frac{\omega_{n-1}}{qAM^{n-1}}\right)^{1/(n-1)}<
\left(\frac{\omega_{n-1}}{4A}\right)^{1/(n-1)},$ так что по лемме
\ref{lem1} точка $b$ не может быть существенно особой для $f.$
Предположим, что $b=0$ является для отображения $f$ полюсом. Тогда
рассмотрим композицию отображений $h=g\circ f,$ где
$g(x)=\frac{x}{|x|^2}$ -- инверсия относительно единичной сферы
${\Bbb S}^{n-1}.$ Заметим, что $h$ -- открытое дискретное кольцевое
$Q$-отображение в нуле и $h(0)=0.$ Кроме того, в некоторой
окрестности нуля отображение $h$ (по построению) является
ограниченным. В таком случае, найдутся $\varepsilon_0>0$ и $R>0,$
такие, что $|h(x)|\le R$ при $|x|<\varepsilon_0.$ Следовательно,
возможно применение предложения \ref{pr3}. По неравенству
(\ref{eq15}), $|h(x)|=\frac{1}{|f(x)|}\le
\frac{\alpha_n(1+R^2)}{\delta}\exp\{-\beta_n I\left(|x|,
\varepsilon_0\right)\}.$ Отсюда следует, что
$$|f(x)|\cdot \exp\{-\beta_n I\left(|x|,
\varepsilon_0\right)\}\ge \frac{\delta}{\alpha_n(1+R^2)}\,.$$
Однако, последнее соотношение противоречит (\ref{eq16}). Полученное
противоречие доказывает, что точка $b=0$ является устранимой для
отображения $f.$ $\Box$
\end{proof}

\medskip
{\bf 3. Основные следствия.} Прежде всего, проведём {\it
доказательство теоремы \ref{th1A}.} В лемме \ref{lem1} полагаем
$\psi(t)=1/tq_b^{1/(n-1)}(t),$ $\varphi(t)=
\exp\left\{\int\limits_{t}^{\varepsilon_0}\frac{dr}{rq_{b}^{1/(n-1)}(r)}\right\}.$
Отметим, что по условию теоремы, при почти всех $r\in (0,
\varepsilon_0)$ имеем $q_b(r)<\infty,$ откуда вытекает строгое
убывание функции $\varphi.$ Кроме того, по теореме Фубини имеем
$$\int\limits_{\varepsilon<|x-b|<\varepsilon_0}Q(x)\cdot\psi^n(|x-b|) dm(x)=$$
$$=\int\limits_{\varepsilon}^{\varepsilon_0}\int\limits_{S(b, r)}
Q(x)\cdot\psi^n(|x-b|)\,dS\,dr =
\omega_{n-1}\cdot\int\limits_{\varepsilon}^{\varepsilon_0}r^{n-1}\psi^n(r)
q_b(r)dr =$$ $$= \omega_{n-1}\cdot I(\varepsilon, \varepsilon_0)=
\omega_{n-1}\cdot\log\varphi(\varepsilon)\,.$$
Отсюда, в частности, следует, что выполнено условие (\ref{eq4A}) при
$M=1.$ Оставшаяся часть утверждения следует из леммы \ref{lem1} и
следствия \ref{cor1}. $\Box$

Полагая в лемме \ref{lem1} в качестве функции $\psi(t)=1/t,$ а в
качестве $\varphi(t)=t^{-q},$ $q>0,$ получаем следующее утверждение.

\medskip
\begin{theorem}\label{th6}
{\sl Пусть $Q:D\rightarrow\,(0\,,\infty]$ -- измеримая по Лебегу
функция, $b\in D$ и $f:D\setminus\{b\}\rightarrow {\Bbb R}^n$ --
открытое дискретное кольцевое $Q$-отображение в точке $b.$
Зафиксируем $\varepsilon_0>0,$ такое что $\varepsilon_0<{\rm
dist}\,(b,
\partial D).$ Пусть, кроме того,
существует число $\varepsilon_0\in(0,1)$ такое, что при
$\varepsilon\rightarrow 0$
$$
\int\limits_{\varepsilon<|x-b|<\varepsilon_0}\frac{Q(x)}{|x-b|^n}\,
dm(x)\le C_1\cdot\log\frac{1}{\varepsilon}\,,$$
где $C_1=C_1(b)$ -- некоторая положительная постоянная.
Тогда существует некоторое число $q>0$ и постоянная $C>0$ такие, что
условие
$$|f(x)|\le C |x-b|^{-q}\,,\qquad \forall x\in
B(b, \varepsilon_0)\setminus\{b\}$$
влечёт наличие предела (конечного либо бесконечного) отображения $f$
в точке $b.$
Кроме того, существует постоянная $p>0,$ такая что оценка вида
$$\lim\limits_{x\rightarrow b}|f(x)|\cdot|x-b|^{p}=0
$$
влечёт, что точка $b$ является устранимой для отображения $f.$}
\end{theorem}

\medskip
\begin{proof}
Заметим, что $I(\varepsilon,
\varepsilon_0)=\log\frac{\varepsilon_0}{\varepsilon},$ где, как и
прежде, $I(\varepsilon, \varepsilon_0)$ задаётся соотношением вида
(\ref{eq11}). Всё остальное непосредственно вытекает из леммы
\ref{lem1} и следствия \ref{cor1}.
\end{proof}$\Box$

\medskip
\begin{theorem}\label{th1}{\sl Пусть $Q:D\rightarrow\,(0\,,\infty]$ -- измеримая по Лебегу функция,
$b\in D$ и $f:D\setminus\{b\}\rightarrow {\Bbb R}^n$ -- открытое
дискретное кольцевое $Q$-отображение в точке $b.$ Зафиксируем
$\varepsilon_0>0,$ такое что $\varepsilon_0<{\rm dist}\,(b,
\partial D).$ Предположим, что при почти
всех $x\in D$ и при некотором $K>0$
\begin{equation}\label{eq27}
q_{b}(r)\le K\cdot \left(\log\frac{1}{r}\right)^{\alpha}
\end{equation}
при $r\rightarrow 0.$ Тогда найдётся число $p>0$ такое, что
неравенство (\ref{eq3}), выполненное  при
$\varphi(t)=\exp\log^{\frac{n-\alpha-1}{n-1}}(1/t)$ и некотором
$\alpha\in (0, n-1)$ (где $p>0$ и $C>0$ -- фиксированные постоянные)
влечёт, что точка $b$ является для отображения $f$ либо полюсом,
либо устранимой особой точкой. }\end{theorem}

\begin{proof} Предположим, $\varepsilon_0>0,$
$\varepsilon_0\in \left(0,\, \min\,\{1, {\rm dist\,}(b,
\partial D)\}\right).$ Полагаем $\psi(t)=1/t,$ $I(\varepsilon, \varepsilon_0):
=\int\limits_{\varepsilon}^{\varepsilon_0}\psi(t)dt.$ Тогда при
некоторой постоянной $C_1>0$ и $\varepsilon\rightarrow 0,$ учитывая
условие (\ref{eq27}),  имеем
$$\frac{1}{I^n(\varepsilon,
\varepsilon_0)}\int\limits_{\varepsilon<|x-b|<\varepsilon_0}\frac{Q(x)dm(x)}{|x-b|^n}\le
C_1\cdot \frac{1}{\log^{n-\alpha-1}(1/\varepsilon)}=
C_1\log^{1-n}\varphi(\varepsilon)\,,$$ откуда вытекает выполнение
условий (\ref{eq4})--(\ref{eq11}) леммы \ref{lem1}. Применяя эту
лемму, заключаем, что точка $b$ не может быть существенно особой для
отображения $f.$
\end{proof}$\Box$

\medskip
\begin{corollary}\label{cor2}{\sl Пусть $Q:D\rightarrow\,(0\,,\infty]$ -- измеримая по Лебегу функция,
$b\in D$ и $f:D\setminus\{b\}\rightarrow {\Bbb R}^n$ -- открытое
дискретное кольцевое $Q$-отображение в точке $b.$ Зафиксируем
$\varepsilon_0>0,$ такое что $\varepsilon_0<{\rm dist}\,(b,
\partial D).$ Тогда найдётся
положительная постоянная $\gamma>0,$ $\gamma<1,$ со следующим
свойством. Если отображение $f$ удовлетворяет более сильному, чем
(\ref{eq3}) (при
$\varphi(t)=\exp\log^{\frac{n-\alpha-1}{n-1}}(1/t)$), условию
\begin{equation}\label{eq29}
\lim\limits_{x\rightarrow b}|f(x)|\cdot
\exp\left\{-\gamma\log^{\frac{n-1-\alpha}{n-1}}(\varepsilon_0/|x-b|)\right\}=0\,,
\end{equation} то $f$ имеет устранимую
особенность в точке $b.$}
\end{corollary}

\begin{proof}
По теореме \ref{th1} точка $b$ не может быть существенно особой для
отображения $f.$ Осталось показать, что при выполнении более сильных
условий (\ref{eq27})--(\ref{eq29}) точка $b$ для отображения $f$
также не может быть полюсом. Заметим, что при достаточно большом
$M_1>0,$ некотором $\varepsilon_1\in (0, \varepsilon_0)$ и всех
$\varepsilon\in (0, \varepsilon_1)$
\begin{equation}\label{eq30}
\int\limits_{\varepsilon<|x-b|<\varepsilon_0}\frac{Q(x)dm(x)}{|x-b|^n}\le
M_1\cdot I^{\alpha+1} (\varepsilon, \varepsilon_0)\,,
\end{equation}
где $\psi(t)=1/t$ и $I(\varepsilon,
\varepsilon_0)=\int\limits_{\varepsilon}^{\varepsilon_0}\psi(t)dt.$
Предположим, что точка $b$ является для отображения $f$ полюсом.
Рассмотрим композицию отображений $h=g\circ f,$ где
$g(x)=\frac{x-b}{|x-b|^2}$ -- инверсия относительно единичной сферы
$S(b, 1).$ Заметим, что отображение $h$ является кольцевым
$Q$-отображением в нуле, и $h(0)=0.$ Кроме того, в некоторой
окрестности нуля отображение $h,$ по построению, является
ограниченным. В таком случае, найдутся $\varepsilon_2>0,$
$\varepsilon_2<\varepsilon_1,$ и $R>0,$ такие, что $|h(x)|\le R$ при
$|x-b|<\varepsilon_2.$ Согласно \cite[лемма~5]{S$_3$}, учитывая
соотношение (\ref{eq30}), получаем
$$|h(x)|=\frac{1}{|f(x)|}\le
\frac{\alpha_n(1+R^2)}{\delta}\exp\left\{-\gamma \cdot
\log^{\frac{n-1-\alpha}{n-1}}\frac{\varepsilon_0}{|x-b|}\right\}\,,$$
где постоянные $\alpha_n$ и
$\gamma=\left(\frac{\omega_{n-1}}{M_1}\right)^{1/(n-1)}$ зависят
только от $n,$ а $\delta$ -- от $R.$ Не ограничивая общности, можно
считать, что $\gamma<1.$ Тогда
$$|f(x)|\cdot \exp\left\{-\gamma \cdot
\log^{\frac{n-1-\alpha}{n-1}}\frac{\varepsilon_0}{|x-b|}\right\}\ge
\frac{\delta}{\alpha_n(1+R^2)}\,.$$ Однако, если при этом выполнено
условие (\ref{eq29}), мы приходим к противоречию, что и доказывает
следствие.
\end{proof}$\Box$

\medskip
Пусть $\Phi:[0,\infty ]\rightarrow [0,\infty ]$ -- неубывающая
функция. Тогда {\it обратная функция} $\Phi^{-1}$ может быть
корректно определена следующим образом:
\begin{equation}\label{eq5.5CC} \Phi^{-1}(\tau)\ =\
\inf\limits_{\Phi(t)\ge \tau} t\,.
\end{equation} Как обычно, $\inf$ в (\ref{eq5.5CC}) равен $\infty,$
если множество $t\in[0,\infty ],$ таких что $\Phi(t)\ge \tau,$
пусто. Заметим, что функция $\Phi^{-1}$ также является неубывающей.

\medskip
\begin{theorem}\label{th7}
{\sl Пусть $Q:D\rightarrow\,(0\,,\infty]$ -- измеримая по Лебегу
функция, $b\in D$ и $f:D\setminus\{b\}\rightarrow {\Bbb R}^n$ --
открытое дискретное кольцевое $Q$-отображение в точке $b.$
Зафиксируем $\varepsilon_0>0,$ такое что $\varepsilon_0<{\rm
dist}\,(b,
\partial D).$
Кроме того, предположим, что существуют число $M>0,$ неубывающая
выпуклая функция $\Phi:[0, \infty]\rightarrow [0, \infty]$ и
окрестность $U$ точки $b,$ такие что
\begin{equation}\label{eq31}
\int\limits_U\Phi\left(Q(x)\right)\frac{dm(x)}{\left(1+|x|^2\right)^n}\
\le\ M
\end{equation}
и
\begin{equation}\label{eq32}
\int\limits_{\delta_0}^{\infty}
\frac{d\tau}{\tau\left[\Phi^{-1}(\tau)\right]^{\frac{1}{n-1}}}=
\infty
\end{equation} при некотором $\delta_0>\Phi(0).$
Тогда найдётся число $p>0,$ такое что неравенство (\ref{eq2}),
выполненное с некоторой постоянной $C>0,$ влечёт наличие предела
отображения $f$  в точке $b.$ }
\end{theorem}

\medskip
\begin{proof}
Из условий (\ref{eq31})--(\ref{eq32}) и \cite[теорема~3.1]{RSev}
следует расходимость интеграла вида
$\int\limits_{0}^{\varepsilon_0}\frac{dt}{tq_{b}^{1/(n-1)}(t)}=\infty.$
Всё остальное следует из теоремы \ref{th1A}.
\end{proof} $\Box$

\medskip Следуя \cite{IR}, будем говорить, что локально интегрируемая
функция ${\varphi}:D\rightarrow{\Bbb R}$ имеет {\it конечное среднее
колебание} в точке $x_0\in D$, пишем $\varphi\in FMO(x_0),$ если
%
%
%
%\begin{equation}\label{eq29*!}
%
$$\limsup\limits_{\varepsilon\rightarrow
0}\frac{1}{\Omega_n\varepsilon^n}\int\limits_{B( x_0,\,\varepsilon)}
|{\varphi}(x)-\overline{{\varphi}}_{\varepsilon}|\,
dm(x)<\infty\,,$$
%
%\end{equation}
%
где
$\overline{{\varphi}}_{\varepsilon}=\frac{1}
{\Omega_n\varepsilon^n}\int\limits_{B(x_0,\,\varepsilon)}
{\varphi}(x)\, dm(x).$
\medskip
Заметим, что, как известно, $\Omega_n\varepsilon^n=m(B(x_0,
\varepsilon)).$ Имеет место следующая

\medskip
\begin{theorem}\label{th2} {\sl Пусть $Q:D\rightarrow\,(0\,,\infty]$ -- измеримая по Лебегу функция,
$b\in D$ и $f:D\setminus\{b\}\rightarrow {\Bbb R}^n$ -- открытое
дискретное кольцевое $Q$-отображение в точке $b.$ Предположим, что
$Q\in FMO(b).$ Тогда найдётся показатель степени $p^{\,\prime},$
зависящий только от функции $Q,$ такой что неравенство
\begin{equation}\label{eq2A}
|f(x)|\le C \left(\log\frac{1}{|x-b|}\right)^{p}\,,
\end{equation}
выполненное при некотором $C>0$ и $p\in (0, p^{\,\prime})$ влечёт
наличие предела отображения $f$ в точке $b.$

Более того, существует $0<p_0<p^{\,\prime},$ такое, что оценка вида
\begin{equation}\label{eq18A}\lim\limits_{x\rightarrow
b}\frac{|f(x)|}{\left(\log\frac{1}{|x-b|}\right)^{p_0}}=0
\end{equation}
влечёт, что $b=0$ является устранимой особой точкой для отображения
$f.$ }
\end{theorem}

\medskip
\begin{proof}
В лемме \ref{lem1} положим $\psi(t):=\frac{1}{t\,\log{\frac1t}},$
т.е., $I(\varepsilon,
\varepsilon_0)=\int\limits_{\varepsilon}^{\varepsilon_0}\psi(t) dt
=\log{\frac{\log{\frac{1}
{\varepsilon}}}{\log{\frac{1}{\varepsilon_0}}}},$ кроме того,
$\varphi(\varepsilon):=\log\frac{1}{\varepsilon}.$ Заметим, что для
функций класса $FMO$ в точке $b$
$$
\int\limits_{\varepsilon<|x|<{e_0}}\frac{Q(x+b)\, dm(x)} {\left(|x|
\log \frac{1}{|x|}\right)^n} = O \left(\log\log
\frac{1}{\varepsilon}\right)
$$
при  $\varepsilon \rightarrow 0 $ и для некоторого $e_0>0,$ $e_0 \le
{\rm dist}\,\left(0,\partial D\right)$ (см.
\cite[следствие~6.3]{MRSY}). Поэтому мы сразу получаем выполнение
условий (\ref{eq4}) и (\ref{eq11}), что доказывает возможность
непрерывного продолжения в точку $b.$ Оставшаяся часть утверждения
вытекает из следствия \ref{cor1}.~$\Box$
\end{proof}

\medskip
\begin{theorem}\label{th1B}{\sl\, Пусть $Q:D\rightarrow\,(0\,,\infty]$ -- измеримая по Лебегу функция,
$b\in D$ и $f:D\setminus\{b\}\rightarrow {\Bbb R}^n$ -- открытое
дискретное кольцевое $Q$-отображение в точке $b.$ Пусть, кроме того,
$q_{b}(r)\le C\cdot \left(\log\frac{1}{r}\right)^{n-1}$ при
$r\rightarrow 0.$ Тогда существует $0<p_0<\infty$ такое, что точка
$b$ является для отображения $f$ либо полюсом, либо устранимой
особой точкой, как только соотношение (\ref{eq2A}) выполнено при
произвольном $0<p<p_0.$

Кроме того, найдётся $0<p_1<\infty$ такое, что если при $0<p<p_1$
имеет место оценка вида (\ref{eq18A}), то точка $b=0$ является
устранимой для отображения $f.$}\end{theorem}

\begin{proof}
Не ограничивая общности рассуждений, можно считать, что $b=0.$
Фиксируем $\varepsilon_0<\min\left\{{\rm dist\,}\left(0,
\partial D\right),\quad 1\right\}.$ Полагаем
$\psi(t)=\frac{1}{t\,\log{\frac{1}{t}}},$
$\varphi(\varepsilon):=\log\frac{1}{\varepsilon}.$  Заметим, что
$$\int\limits_{\varepsilon<|x|<\varepsilon_0}
\frac{Q(x)dm(x)}{\left(|
x|\log{\frac{1}{|x|}}\right)^n}=\int\limits_{\varepsilon}^{\varepsilon_0}
\left(\int\limits_{|x|\,=\,r}\frac{Q(x)dm(x)}{\left(|
x|\log{\frac{1}{|x|}}\right)^n}\, dS\,\right)\,dr\le C\cdot\omega_
{n-1}\cdot I(\varepsilon, \varepsilon_0)\,,$$
где, как и прежде, $I(\varepsilon,
\varepsilon_0):=\int\limits_{\varepsilon}^{\varepsilon_0}\psi(t) dt
= \log{\frac{\log{\frac{1}
{\varepsilon}}}{\log{\frac{1}{\varepsilon_0}}}}.$ Отсюда заключаем,
что при указанной выше функции $\psi$ имеют место условия
(\ref{eq4}) и (\ref{eq11}) леммы \ref{lem1}. Таким образом, первая
часть заключения теоремы \ref{th1B} установлена. Второе утверждение
этой теоремы вытекает из следствия \ref{cor1}.
\end{proof}$\Box$

\medskip
Напоследок сформулируем ещё одно важное утверждение, относящееся к
одному достаточно широкому классу отображений, для которого основные
результаты заметки имеют место. Напомним, что точка $x_0\in D$
называется {\it точкой ветвления отображения $f$}, если ни в какой
окрестности $U$ точки $x_0$ сужение $f|_U$ не является
гомеоморфизмом. Множество точек ветвления отображения $f$ принято
обозначать символом $B_f.$ Поскольку  произвольное открытое
дискретное отображение $f:D\rightarrow {\Bbb R}^n$ класса
$W_{loc}^{1,n}(D),$ для которого $K_I(x,f)\in L_{loc}^1$ и
$m(B_f)=0,$ является кольцевым $Q$-отображением в произвольной точке
заданной области, см. \cite[теорема~1]{S$_4$} и \cite[следствие~3.14
и предложение~4.3]{MRSY$_1$}), из теорем \ref{th1A}--\ref{th7}
вытекает следующее утверждение.

\medskip
\begin{corollary}\label{cor3}{\sl
Каждое из заключения теорем \ref{th1A}--\ref{th7}, а также следствий
\ref{cor1}--\ref{cor2}, справедливо для произвольных открытых
дискретных отображений $f:D\setminus\{b\}\rightarrow {\Bbb R}^n$
класса $W_{loc}^{1,n}(D\setminus\{b\}),$ таких что $K_I(x,f)\in
L_{loc}^1$ и $m(B_f)=0.$}
\end{corollary}

КОНТАКТНАЯ ИНФОРМАЦИЯ

\medskip
\noindent{{\bf Дарья Сергеевна Доля} \\
Институт математики НАН Украины \\
ул. Терещенковская, 3 \\
г.~Киев-4, Украина, 01 601 \\ тел. +38 066 889 46 50 (моб.), e-mail:
dasha.dolja@mail.ru}

\medskip
\noindent{{\bf Евгений Александрович Севостьянов} \\
Житомирский государственный университет им.\ И.~Франко\\
кафедра математического анализа, ул. Большая Бердичевская, 40 \\
г.~Житомир, Украина, 10 008 \\ тел. +38 066 959 50 34 (моб.),
e-mail: esevostyanov2009@mail.ru}

\end{document}